# The exception of generality quantifier in the formulas of the elementary theory in a free group.

**G.S. Makanin**

The method of exception of generality quantifier in the formulas of the elementary theory in a free group is offered. The solvability of the elementary theory of a free group is proved.

Let **G** be a free group with free generators $\mathbf{a_1,...,a_\omega}$.
Let $\mathbf{x_1,...,x_k}$ be *variables*, whose values are elements of the free group **G.** *The equation* in a free group is named an equality, in the form

$$\mathbf{W(x_1,...,x_k, a_1,...,a_\omega) = 1} , \qquad \mathbf{W}$$

where **W** is a word, composed by letters $\mathbf{x_1,...,x_k, a_1,...,a_\omega}$ and their inverses.
The generators $\mathbf{a_1,...,a_\omega}$ are called the *coefficients* of the equation.

A list of words $\mathbf{X_1 ,..., X_k \in \langle a_1,...,a_\omega \rangle}$ is called a solution of the equation **W**, if $\mathbf{W(X_1,...,X_k, a_1,...,a_\omega) = 1}$ in the free group **G**.

In **[1],[2]** there is an algorithm constructed, recognizing a decidability of equations in a free group. The decidability of the universal and positive theories is established.

The elementary theory of the free group G is a universe of true closed formulas of restricted predicate calculus of a group signature
$(=,\cdot,^{-1},\mathbf{a_1} ,...,\mathbf{a_\omega})$. The problem of decidability of a theory of a free group was set by A.Tarsky in 1957-1959 years.

It is known that every closed formula of the restricted predicate calculus with equality in G can be represented in the form

$$\mathbf{Q_1x_1...Q_kx_k( \bigvee_{i=1}^{q}(\&_{j=1}^{r_i} W_{i,j} (x_1,..., x_k, a_1,..., a_\omega) = \neq\ 1))}$$

where each $\mathbf{Q_i}$ is a universal or existential quantifier; $\mathbf{x_1,...,x_k}$ are variables; $\mathbf{a_1,..., a_\omega}$ are generators of G-group; each atomic formula is either an equation $\mathbf{W_{i,j}(x_1,..., x_k, a_1,..., a_\omega) = 1,}$ or an inequality $\mathbf{W_{i,j}(x_1,..., x_k, a_1,..., a_\omega) \neq 1}$

*The formula of the form*

$$\mathbf{Q_1x_1...Q_kx_k\exists z ( W(x_1,...,x_k,z,a_1,..., a_\omega) = 1\ \&z\neq 1)} , \qquad \mathbf{R_1}$$

*where* $\mathbf{k\geq 1}$ , $\mathbf{Q_i}$ *are quantifiers, is called reduced.*

In **[3]** it is proved, that its equivalent reduced formula can be constructed for each closed formula of restricted predicate calculus with an equality in G.

**The way of reduction.**

The way of reduction of the equation $\mathbf{W(x_1...,x_{k-1},z,a_1, ...,a_\omega) = 1}$ in the free group **G** divides the word $\mathbf{W(x_1...,x_{k-1},z,a_1, ...,a_\omega)}$ into parts (variables can be divided) and indicates what parts should be reduced among themselves to reduce the word **W( … )** completely.
The equation will have the form of $\mathbf{g_1,...,g_m = 1}$ where $\mathbf{g_1,...,g_m}$ are the designations of the parts of the equation.

*The assignment* $\zeta(g_i)$ of the letter $g_i$ of the equation $g_1...g_m= 1$ is called a word, which consists of dictionary variables $\upsilon_i^{\varepsilon}$ and coefficients $\alpha_j^{\delta}$. The letter $a_p^{\delta}$ of the equation always has as an assignment the coefficient $\alpha_p^{\delta}$.

*The arc of reduction* of the equation **W** with the way of reduction $\Sigma$ is called a tape, which connects *from above* (taking entirely) the parts of the equation. The arc can be linear.

Each pair of the parts of the equation can be connected by no more than one general arch.

The arches of an equation do not cross among themselves.

The way of reduction is a pair of connected among themselves lines

$$\begin{matrix} \zeta(g_1) & \zeta(g_2)... & \zeta(g_m) \\ g_1 & g_2 & g_m \end{matrix}$$

where there are parts of the equation $g_1...g_m= 1$ , in the lower line, and there are assignments $\zeta(g_1),\zeta(g_2),\ldots,\zeta(g_m)$ of the letters $g_1,g_2,...,g_m$ in the higher line, and at the same time any $\zeta(g_i)$ is written over the occurrence $g_i$.

The product of all the assignments $\zeta(g_1)\ldots\zeta(g_m)$ is equal to **1** in the free group $G'\langle\upsilon_1,...,\upsilon_r,\alpha_1,...,\alpha_\omega\rangle$.

Solution of the equation **W** presupposes for each $i=1,...,m$ the equality $\zeta(g_i)=g_i$ in the group $G'$.

If the assignment $\zeta(g_i)$ of the variable $g_i$ is empty, then $g=1$ in $G'$.

If an arch (tape, linear) covers from above the subword $S(x_1...,x_k,a_1, ...,a_\omega)$ , then $S(x_1...,x_k,a_1, ...,a_\omega)=1$ in G.

**Lemma 1**. For each reduced formula $R_1$ its equivalent disjunction of reduced formulas with the given way of reduction can be constructed.

**Proof.** Let's consider the reduced formula $R_1$

$$Q_1x_1\ldots Q_{k-1}x_{k-1}\exists z\ (W(x_1,...,x_{k-1}, z, a_1,...,a_\omega)=1\&z\neq1)$$

Each pair of parts of the equation $W(x_1,...,x_k, z, a_1,...,a_\omega)=1$ can be connected by no more than one arch. The first part $g_1$ of the equation can be connected with all other parts by no more, than **d-1** arches.

The second part $g_2$ of the equation can be connected with all other parts also by no more, than **d-1** arches, etc. Thus, there are no more than $(d-1)^d$ different ways of reduction of the equation, which length is equal to **d**.

Let's construct all kinds of reduction of the reduced formula $R_1$: $\Sigma_1$ ,..., $\Sigma_s$. And finally, we will construct a disjunction of reduced formulas with given way of reduction $R_1\&\Sigma_1 V...VR_1\&\Sigma_s$.

Let $\mathbf{R_1}$ be a reduced formula with the way of reduction $\mathbf{\Sigma_1}$

$$\mathbf{Q_1x_1 \ldots Q_{k-1}x_{k-1} \exists z\, ((W(x_1,\ldots,x_{k-1},a_1,\ldots,a_\omega)=1 \& \Sigma_1 \& z \neq 1)} \qquad \mathbf{R_1 \& \Sigma_1}$$

where $\mathbf{Q_1,\ldots,Q_{k-1}}$ are quantifiers $\forall$ and $\exists$, where (*here and from now on*) all $\forall$-variables $\mathbf{x_{\lambda_1},\ldots, x_{\lambda_\rho}}$ have another definition (respectively) $\mathbf{y_{\lambda_1},\ldots, y_{\lambda_\rho}}$; all $\exists$-variables $\mathbf{x_{\mu_1},\ldots, x_{\mu_\sigma}}$ have another definition $\mathbf{z_{\mu_1},\ldots, z_{\mu_\sigma}}$, $\mathbf{\Sigma_1}$ is a way of reduction of the formula's equation $\mathbf{R_1\&\Sigma_1}$.

It's obvious that if $\forall$-variable $\mathbf{y_i^\varepsilon}$ of the formula $\mathbf{R_1\&\Sigma_1}$ connects by its assignment with an assignment of a word in an alphabet $\langle \mathbf{a_1,\ldots,a_\omega} \rangle$, than the formula is false.

Let non self-inverse $\forall$-variables $\mathbf{y_i^\varepsilon}$ and $\mathbf{y_j^\delta}$ be connected by their assignments. Because a solution of one assignment can be any word in an alphabet $\mathbf{(a_1,\ldots,a_\omega)}$ , and a solution of another assignment can be any word in an alphabet $\mathbf{(a_1,\ldots,a_\omega)}$, then this connection cannot exist.

If self-inverse variables $\mathbf{y_i}$ and $\mathbf{y_i^{-1}}$ are connected by their assignments, then they are connected entirely and will be called a *superfluous pair*.

It is clear that if all superfluous pairs are removed (in any order) from the formula with superfluous pairs, then we will obtain its equivalent formula without superfluous pairs.

*A reduced formula without superfluous pairs is called a prepared formula.*

Let $\mathbf{R_2}$ be a prepared formula with the way of reduction $\mathbf{\Sigma_2}$

$$\mathbf{Q_1x_1 \ldots Q_{k-1}x_{k-1} \exists x_k((W(y_1,\ldots,y_\rho,z_1,\ldots,z_\sigma,a_1,\ldots,a_\omega)=1 \& \Sigma_2 \& x_k \neq 1)} \qquad \mathbf{R_2 \& \Sigma_2}$$

where each variable $\mathbf{y}$ is connected *only* with variables $\mathbf{z}$.

## Definition of the fragment

*A platform* of the prepared formula $\mathbf{R_2\&\Sigma_2}$ is called a subword $\mathbf{P}$ of the word $\mathbf{W(y_1,\ldots,y_\rho,z_1,\ldots,z_\sigma,a_1,\ldots,a_\omega)}$ in the form

$$\mathbf{z_{\beta_1}^{\delta_1} K_1 z_{\beta_2}^{\delta_2} K_2 \ \ldots\ z_{\beta_q}^{\delta_q} K_q \cdot y_j \cdot K_{q+1}\, z_{\beta_{q+1}}^{\delta_{q+1}} \ldots K_{s-1}\, z_{\beta_{s-1}}^{\delta_{s-1}}\, K_s\, z_{\beta_s}^{\delta_s}} \qquad \mathbf{P}$$

where

$\mathbf{\delta_2 = \pm 1}$, $\mathbf{y_j}$ is the only $\mathbf{y}$-variable of the platform $\mathbf{P}$ ;

$\mathbf{y_j}$ is connected by arches with all $\mathbf{z_{\beta_i}^{\delta_i}}$ $\mathbf{(i=1,\ldots,s)}$ of the platform $\mathbf{P}$ and only with them;

arches of the variable $\mathbf{y_j}$ isolate in variables $\mathbf{z_{\beta_i}^{\delta_i}}$ variables $\mathbf{c_i}$ and equalities $\mathbf{z_{\beta_i}^{\delta_i} = r_i\, c_i\, t_i}$ $\mathbf{(i=1,\ldots,s)}$;

the variable $\mathbf{y_j}$ is reduced entirely with the foundations of its arches, the reduction $\mathbf{y_j}$ defines the equality

$$\mathbf{y_j = c_q^{-1}\, c_{q-1}^{-1} \ldots c_2^{-1}\, c_1^{-1}\, c_s^{-1}\, c_{s-1}^{-1} \ldots c_{q+2}^{-1}\, c_{q+1}^{-1}} .$$

Let's call part $\mathbf{r_1}$ in the variable $\mathbf{z_{\beta_1}^{\delta_1} = r_1 c_1 t_1}$ *detached*. Let's call part $\mathbf{t_s}$ in the variable $\mathbf{z_{\beta_s}^{\delta_s} = r_s c_s t_s}$ *detached.*

The fragment $\mathbf{\Phi}$ of the word $\mathbf{W}$ $(\mathbf{y_1,\dots,y_\rho,z_1,\dots,z_\sigma,a_1,\dots,a_\omega})$ is called the platform $\mathbf{P}$ without separate parts: without letter $\mathbf{r_1}$ in $\mathbf{r_1c_1t_1}$ and without letter $\mathbf{t_s}$ in $\mathbf{r_sc_st_s}$.

A detailed *definition of the fragment* $\mathbf{\Phi}$ with the variable $\mathbf{y_j}$:

$$\mathbf{c_1t_1K_1r_2c_2\,t_2\,K_2\dots r_{q-1}c_{q-1}t_{q-1}\,K_{q-1}r_qc_q\,t_q\,K_q\cdot y_j\cdot K_{q+1}\,r_{q+1}c_{q+1}\,t_{q+1}\,K_{q+2}\,r_{q+2}\,c_{q+2}\,t_{q+2}\dots K_{s-1}\,r_{s-1}\,c_{s-1}\,t_{s-1}\,K_s\,r_s\,c_s} \qquad \mathbf{\Phi}$$

where

$\mathbf{c_1c_2\dots c_{q-1}\,c_q\,y_jc_{q+1}c_{q+2}\dots c_{s-1}c_s=1}$ in $\mathbf{G}$;

$\mathbf{y_i}$ is connected with all $\mathbf{c_i}$ and absorb them.

The *fully* definition of the fragment $\mathbf{\Phi}$ will be obtained if we substitute in a detailed fragment $\mathbf{\Phi}$ the variable $\mathbf{y_i}$ for $\mathbf{c_q^{-1}\,c_{q-1}^{-1}\dots c_2^{-1}\,c_1^{-1}\,c_s^{-1}\,c_{s-1}^{-1}\dots c_{q+2}^{-1}\,c_{q+1}^{-1}}$.

The arches of the fully fragment $\mathbf{\Phi}$ that connect $\mathbf{y_i}$ with $\mathbf{c_1,\dots,c_q}$, $\mathbf{c_{q+1},\dots,c_s}$ isolate the conjunction of equalities in the group $\mathbf{G}$

$$
\begin{aligned}
&\mathbf{H_1= t_1\,K_1r_2=1} \quad \& \\
&\mathbf{H_2= t_2\,K_2r_3=1} \quad \& \\
&\dots\dots\dots\dots \\
&\mathbf{H_{q-1}=t_{q-1}K_{q-1}r_q=1} \quad \& \\
&\mathbf{H_q\ =t_q\,K_q=1} \quad \& \\
&\mathbf{H_{q+1}=K_{q+1}r_{q+1}=1}\,\& \\
&\mathbf{H_{q+2}=t_{q+1}\,K_{q+2}r_{q+1}=1}\,\& \\
&\dots\dots\dots\dots\dots\dots \\
&\mathbf{H_{s-1}=t_{s-2}\,K_{s-1}r_{s-1}=1}\,\& \\
&\mathbf{H_s=t_{s-1}\,K_sr_s=1}
\end{aligned}
$$

**<u>Lemma 2</u>**. The prepared formula $\mathbf{R_2\&\Sigma_2}$ contains a platform.

**Proof**. A prepared formula doesn't contain superfluous pairs (the united variables $(\mathbf{y_i}, \mathbf{y_i^{-1}})$).

If $\mathbf{R_2\&\Sigma_2}$ doesn't contain $\mathbf{y}$-variables, than $\mathbf{R_2\&\Sigma_2}$ is a universal formula.

To choose a platform with the only $\mathbf{y}$-variable in the word $\mathbf{W(y_1,\dots,y_\rho,z_1,\dots,z_\sigma,a_1,\dots,a_\omega)}$, let's choose any occurrence of the $\mathbf{y}$-variable, let $\mathbf{y_1}$. Let's isolate all connections of the chosen $\mathbf{y_1}$ with $\mathbf{z}$ -variables. We'll get a word, let this one

$$\mathbf{z_1\,K_1z_2\,K_2\,\dots\,z_r\,K_r\,\,y_1\,K_{r+1}\,z_{r+1}\dots K_m\,z_m} \qquad \mathbf{P}$$

where $\mathbf{y_1}$ must be reduced completely connecting with all $\mathbf{z_{r+1},\dots,z_m}$.

If no one from $\mathbf{K_i}$ contains $\mathbf{y}$-variable, then a required platform is constructed. Let any $\mathbf{K_i}$ from $\mathbf{P_1}$ contain $\mathbf{y_2}$. Then $\mathbf{y_2}$ must be reduced entirely inside $\mathbf{K_i}$. Let's isolate all the connections of the variable $\mathbf{y_2}$ with $\mathbf{z}$-variables etc. Continuing this process, we will get a platform with the only $\mathbf{y}$-variable. We will define it like $\mathbf{P_1}$.

So $\mathbf{W(y_1,\dots,y_\rho,\,z_1,\dots,z_\sigma,a_1,\dots,a_\omega)=W_1(y_1,\dots,\,z_1,\dots,a_1,\dots,a_\omega)P_1W_2(y_1,\dots,\,z_1,\dots,a_1,\dots,a_\omega)}$, where $\mathbf{P_1}$ is a platform with the only $\mathbf{y}$-variable.

**Lemma 3**. The prepared formula $\mathbf{R_2\&\Sigma_2}$ contains a fragment.

**Proof.** By definition there are no superfluous pairs in the word $\mathbf{W(y_1,...,y_\rho, z_1,...,z_\sigma,a_1,...,a_\omega)}$ of a prepared formula. Therefore, **y-**variables can be connected only with **z-**variables

Let's isolate in the word $\mathbf{W(y_1,...,y_\rho,z_1,...,z_\sigma,a_1,...,a_\omega)}$ all **y-**variables and all arches of the **y-**variables.

Let's isolate in the word $\mathbf{W(y_1,...,y_\rho,z_1,...,z_\sigma,a_1,...,a_\omega)}$ some platform with **y-**variable (let it be the variable $\mathbf{y_i}$).

$$z_{\beta_1}^{\delta_1} K_1 z_{\beta_2}^{\delta_2} K_2 \ldots z_{\beta_{q-1}}^{\delta_{q-1}} K_{q-1} z_{\beta_q}^{\delta_q} K_q \cdot y_i \cdot K_{q+1} z_{\beta_{q+1}}^{\delta_{q+1}} K_{q+2} z_{\beta_{q+2}}^{\delta_{q+2}} \ldots K_{s-1} z_{\beta_{s-1}}^{\delta_{s-1}} K_s z_{\beta_s}^{\delta_s} \qquad \mathbf{P}$$

where $z_{\beta_i}^{\delta_i} = r_i\, c_i\, t_i$ $\quad(i=1,...,s)$; $\mathbf{y_i}$ is connected with $\mathbf{c_1,..., c_s}$ and only with them**.**

From the equality **P** we will get the equality

$$y_i = c_q^{-1} c_{q-1}^{-1} \ldots c_2^{-1} c_1^{-1} c_s^{-1} c_{s-1}^{-1} \ldots c_{q+2}^{-1} c_{q+1}^{-1},$$

Let's substitute this equality into **P**. We will get that for all $\alpha=1,\ldots,s$

$c_\alpha$ connects $c_\alpha^{-1}$ from $y_i$ .

and that subword $\mathbf{K_i}$ doesn't contain occurrences of the **y**-variables $(i=1,...,s)$

In the first variable $z_{\beta_1}^{\delta_1} = r_1\, c_1\, t_1$ we will divide part $\mathbf{r_1}$, and in the last variable $z_{\beta_s}^{\delta_s} = r_s\, c_s\, t_s$ we will divide $\mathbf{t_s}$ . We will obtain

$$\Phi = c_1\, t_1 K_1\, r_2 c_2\, t_2\, K_2 \ldots r_q c_q\, t_q\, K_q\, y_i K_{q+1}\, r_{q+1} c_{q+1} t_{q+1} \ldots K_{s-1}\, r_{s-1} c_{s-1}\, t_{s-1}\, K_s\, r_s c_s$$

Here $\mathbf{\Phi}$ is a required fragment. As because of the reduction $\mathbf{y_i}$ with variables $\mathbf{c_\alpha}$ of the fragment $\mathbf{\Phi}$ is entirely covered by arches, then

$\mathbf{\Phi = 1}$ in **G**.

The subword $\mathbf{c_1...c_q\ y_i\ c_{q+1}...c_s}$ we will call a *core* of the fragment $\mathbf{\Phi}$.

The core of the fragment $\mathbf{\Phi}$ consist of mutually connected variables. It can be divided from the other part of the fragment.

We will enter additional designations and equalities

$\mathbf{C_1 = c_1c_2...c_qy_1\ c_{q+1}...c_{s-1}c_s}$, where $\mathbf{C_1=1}$ in **G**

.....................................................................

$\mathbf{C_r = c_pc_{p+1}...c_my_r\ c_{m+1}...c_{t-1}c_t}$ where $\mathbf{C_r=1}$ in **G**

The word $\mathbf{H_1H_2...H_{q-1}H_qH_{q+1}H_{q+2}\ ...H_{s-1}H_s}$ we will define with $\mathbf{D(z_1,...,z_\sigma,a_1,...,a_\omega)}$.

Using $\mathbf{H_i=1\ (i=1,...,s)}$, we will transfer will transfer one after another every $\mathbf{H_i}$ from $\mathbf{\Phi}$ to the end of the word. We will obtain $\mathbf{\Phi = c_1c_2\ ...\ c_{q-1}\ c_q\ y_jc_{q+1}c_{q+2}\ ...\ c_{s-1}c_s\, D(z_1,...,z_\sigma,a_1,...,a_\omega)}$ in the group **G** . It is clear that $\mathbf{D(z_1,...,z_\sigma,a_1,...,a_\omega)=1}$ in **G** .

## **Normal formula**.

A prepared formula in the form of $\mathbf{R_2\&\Sigma_2}$

$$\forall y_1 \ldots \forall y_\rho \exists z_1 \ldots \exists z_\sigma \exists z (\Phi_1 \Phi_2 \ldots \Phi_r\, E(z_1,...,z_\sigma,a_1,...,a_\omega) = 1 \ \&\ \Sigma_2 \ \&\ z \neq 1) \qquad \mathbf{N_1}$$

where $\mathbf{\Phi_1\Phi_2...\Phi_r}$ are fragments, is called a normal and is defined like **N**.

**Lemma 4**. A prepared formula is equivalent to the normal.

**Proof**. Let's consider a prepared formula $R_2\&\Sigma_2$.
Let's isolate in the word $W(y_1,\dots,y_\rho, z_1,\dots,z_\sigma,a_1,\dots,a_\omega)$ the fragment $\Phi_1$ .

Because the fragment $\Phi_1=1$ in $G$ , it can be taken away from $W(y_1,\dots,y_\rho, z_1,\dots,z_\sigma,a_1,\dots,a_\omega)$ and removed to the place *after this word*. We will get the word $W_1(y_1,\dots,y_\rho, z_1,\dots,z_\sigma,a_1,\dots,a_\omega)$ and equality $W(y_1,\dots,y_\rho,z_1,\dots,z_\sigma,a_1,\dots,a_\omega)=W_1(y_1,\dots,y_\rho, z_1,\dots,z_\sigma,a_1,\dots,a_\omega)\Phi_1$ in $G$. Then in the word $W_1(\dots)$ we will isolate the fragment $\Phi_2$, transfer $\Phi_2$ to the place after the word $\Phi_1$ and obtain the equality $W_1(\dots)=W_2(\dots)\Phi_2$, i.e. $W(y_1,\dots,y_\rho, z_1,\dots,z_\sigma,a_1,\dots,a_\omega)=W_2(y_1,\dots,y_\rho, z_1,\dots,z_\sigma,a_1,\dots,a_\omega)\Phi_1\Phi_2$.

Continuing the process of isolating the fragment up to the end, we will get the equality

$$W(y_1,\dots,y_\rho, z_1,\dots,z_\sigma,a_1,\dots,a_\omega)=E(z_1,\dots,z_\sigma,a_1,\dots,a_\omega)\Phi_1\Phi_2\dots\Phi_r \text{ in } G$$

where $E(z_1,\dots,z_\sigma,a_1,\dots,a_\omega)$ is a subword, *that doesn't consist* $y$*-variables*, $E(y_1,\dots,y_\rho, z_1,\dots,z_\sigma,a_1,\dots,a_\omega)=1$ in $G$. Then $W= E (z_1,\dots,z_\sigma,a_1,\dots,a_\omega) C_1D_1\dots C_rD_r$ in $G$.

Because $E(\dots)=1$ in $G$ for every $i=1,\dots,r$ in the word $C_i$ the variable $y_i$ is connected with the list of $c$-letters and is swallowed by them, then in the formula $R_2\&\Sigma_2$ can be placed by the way of

$$\forall y_1\dots\forall y_\rho\exists z_1\dots\exists z_\sigma\exists z\, (E( z_1,\dots,z_\sigma,a_1,\dots,a_\omega)\, C_1D_1\dots C_rD_r=1 \,\&\, \Sigma_2 \,\&\, z\neq 1)$$

So the prepared formula $R_2\&\Sigma_2$ is equivalent to the formula $N_1$ in the form of

$$\forall y_1\dots\forall y_\rho\exists z_1\dots\exists z_\sigma\exists z(\Phi_1\dots\Phi_rE(z_1,\dots,z_\sigma,a_1,\dots,a_\omega) =1\&\Sigma_2\&z\neq 1) \qquad N_1$$

## Matrixes.

(Further so we can manage without negative degrees, the variables $y_1^{-1},\dots,y_\rho^{-1}$ we will define respectively $y_{\rho+1},\dots,y_{\rho+\rho}$ ; the variables $z_1^{-1},\dots,z_\sigma^{-1}$ - $z_{1+\sigma},\dots,z_{\sigma+\sigma}$ ).

Let's examine the formula $N_1$ .

The quantifier-free part

$$\Phi_1 \dots \Phi_rE(z_1,\dots,z_\sigma,a_1,\dots,a_\omega)=1 \qquad M_1$$

Of the formula $N_1$ we will call a *matrix* $M_1$.

1.All $y$-variables of the formula $N_1$ are situated in the matrix $M_1$. Because matrix $M_1$ doesn't consist superfluous pairs of $y$-variables, then every $y$-variable is connected by arches with $c$-variables from $z$-variables and only with them.

Let's isolate from the fragments $\Phi_1 \dots \Phi_r$ of the matrix $M_1$ cores

$$\Phi_1 \dots \Phi_r E = C_1\dots C_rD_1\dots D_rE(z_1,\dots,z_\sigma,a_1,\dots,a_\omega) =1$$

In particular,

$$C_1\dots C_r=1\& D_1\dots D_r E =1$$

We have isolated all occurrences of the variables $z_i$ from $\Phi_1 \dots \Phi_r$ , which are connected by arches with $y$-variables $y_1\dots y_\rho$ .

Let it be $z_{\alpha_1},\dots, z_{\alpha_\mu}$ with isolated $c$-variables $z_{\alpha_i}= r_{\alpha_i}\, c_{\alpha_i}\, t_{\alpha_i}$ $(i=1,\dots,\mu)$, which connect $c$-variables with $y$-variables. These and $c_{\alpha_i}$, and $z_{\alpha_i}$ we will call **1**-*internal variables*. Let $z_{\alpha_{\mu+1}},\dots, z_{\alpha_\nu}$ - (*not* **1**-*internal*) occurrences $z$-variables.

In every $\Phi_i$ from the matrix $M_1$ we will reduce *the core*, i.e. reduce
$C_i \to 1$, $i=1,\ldots,\rho$ .

As a result, we will construct the formula $N_2$, which quantifier-free part we will call the matrix $M_2$.

The length (a number of occurrences of **c**-variables) of the matrix $N_2$ is less than the length of the matrix $N_1$.

2. Let $z_{\beta_1},\ldots,z_{\beta_r}$ be all occurrences of the **z**-variables from $M_2$, each from those is equal to (or inverse) **1**-internal variable from the list $z_{\alpha_1},\ldots,z_{\alpha_\mu}$ (from the matrix $M_1$) but they don't enter into the list of the variables $z_{\alpha_1},\ldots,z_{\alpha_\mu}$.

Each occurrence $z_{\beta_i} \to r_{\beta_i}\, c_{\beta_i}\, t_{\beta_i}$ from the matrix $M_2$, which coincide with some **1**-internal $z_{\alpha_i}$, but is not internal $z_{\alpha_i}$, we will call
**2**-*internal*. Each occurrence $c_{\beta_i}$ from $z_{\beta_i}$ coincides with some
**1**-internal $c_{\alpha_i}$. We will call it **2**-*internal*.

(*Superfluous* couples consisted from **2**-internal **z**-variables we will delete.)
We will notice that not superfluous **2**-internal **c**-variables $c_{\beta_1},\ldots,c_{\beta_r}$ can take any value in an alphabet $\langle a_1,\ldots,a_\omega \rangle$. We *can* define them with $y_{\beta_1},\ldots,y_{\beta_r}$.

When we put quantifiers $\forall$ in a position before variables $y_{\beta_1},\ldots,y_{\beta_r}$ and quantifiers $\exists$ before **c**-variables we will get the formula $N_2'$ and formula $M_2'$ .

In the matrix $M_2'$ we will reduce *the core*.

As a result we will construct the formula $N_3$, which quantifier-free part we will call the matrix $M_3$.

**3**. Let $c_{\gamma_1},\ldots, c_{\gamma_s}$ are all occurrences of **c**-variables from $M_3$, which enter into the list of **z**-variables, equal (or inverse) **2**-internal variables $z_{\beta_1},\ldots, z_{\beta_r}$ from matrix $M_2$, but don't enter into the list of variables $z_{\beta_1},\ldots, z_{\beta_r}$ .

Each occurrence $z_{\gamma_i} \to r_{\gamma_i}\, c_{\gamma_i}\, t_{\gamma_i}$ from matrix $M_3$, that don't coincide with some internal $z_{\beta_i}$, but is not internal $z_\beta$, we will call
**3**-*internal* and i.e. ...

………………………………………………………………………………

**n**. Continuing this process, we will construct a sequence of matrixes

$$M_1, M_2,\ldots, M_n$$

which length decreases every step and on **n** step we will get the matrix $M_n$, in which it is impossible from **z**-variables to get **y**-variables.

The received matrix $M_n$ with quantifiers $\exists\, z_1\ldots\exists z_\mu \exists z$ forms a universal formula. We will define it like $R_n \& \Sigma_n$.

<u>**Lemma 5**</u>. About formula $R_n$ with the way of reduction $\Sigma_n$ we can find out if it is true or false.

**Proof**. Let's examine the formula $R_n \& \Sigma_n$

$$\exists\, z_1\ldots\exists z_\sigma \exists z(W(z_1,\ldots,z_\sigma,a_1,\ldots,a_\omega)=1\ \&\ \Sigma_n \&\ z\neq 1)\ .$$

For every matrix $\mathbf{M_n}$ of the arch with foundations $\mathbf{P_i}$ ,$\mathbf{S_i}$ **i=1,...,w** we will construct a system of equations

$$\mathbf{P_i S_i = 1 \quad (i=1,...,w)} \qquad \mathbf{(1)}$$

Then in the formula $\mathbf{R_n \& \Sigma_n}$ every arch $\mathbf{(P_i ,S_i)}$ destroy. As a result we will contain a conjunction **K** of universal formulas' systems without arches which is equivalent to $\mathbf{R_n \& \Sigma_n}$. Using formulas of integration conjunctions and disjunctions (see **[2], p.738** ), and decidability of universal theory (see), we will find out if the conjunction **K** of formulas' systems true or not.

**Theorem.** The elementary theory of a free group is solvable.

**Proof.** The elementary theory of a free group is solvable. The theorem follows from the lemmas **1-5**.

**Literature.**